\documentclass[12pt, reqno]{amsart}
\usepackage{amsmath, amsthm, amscd, amsfonts, amssymb, graphicx, color}
\usepackage[bookmarksnumbered,plainpages]{hyperref}
\usepackage[all]{xy}

\usepackage{xcolor}
\newtheorem{theorem}{Theorem}[section]
\newtheorem{lemma}[theorem]{Lemma}
\newtheorem{proposition}[theorem]{Proposition}
\newtheorem{corollary}[theorem]{Corollary}
\theoremstyle{definition}

\newtheorem{example}[theorem]{Example}

\theoremstyle{remark}
\newtheorem{remark}[theorem]{Remark}
\numberwithin{equation}{section}
\begin{document}
\setcounter{page}{1}
\title[Some refinements of numerical radius inequalities]{Some refinements of numerical radius inequalities}
\author[Heydarbeygi, Amyari and Khanegir]{Zahra Heydarbeygi, Maryam Amyari$^*$ and Mahnaz Khanehgir}
\address{Department of Mathematics,  Mashhad Branch,
 Islamic Azad University, Mashhad 91735,
Iran.}

\email{zheydarbeygi@yahoo.com}
\email{maryam\_amyari@yahoo.com and amyari@mshdiau.ac.ir}
\email{khanehgir@mshdiau.ac.ir}
\keywords{ bounded linear operator, Hilbert space, norm inequality, numerical radius.\\
*Corresponding author}
\begin{abstract}
In this paper, we give some refinements for the second inequality in $\frac{1}{2}\|A\| \leq w(A) \leq \|A\|$, where $A\in B(H)$.
In particular, if $A$ is hyponormal by refining the Young inequality with the Kantorovich constant $K(\cdot, \cdot)$,
we show that $w(A)\leq \dfrac{1}{\displaystyle {2\inf_{\| x \|=1}}\zeta(x)}\| |A|+|A^{*}|\|\leq \dfrac{1}{2}\| |A|+|A^*|\|$,
where $\zeta(x)=K(\frac{\langle |A|x,x \rangle}{\langle |A^{*}|x,x \rangle},2)^{r},~~~r=\min\{\lambda,1-\lambda\}$ and $0\leq \lambda \leq 1$ .
We also give a reverse for the classical numerical radius power inequality $w(A^{n})\leq w^{n}(A)$
 for any operator $A \in B(H)$ in the case when $n=2$.
\end{abstract}
\maketitle
\section{Introduction}
Suppose that $(H, \langle \cdot,\cdot \rangle)$ is a complex Hilbert space and $B (H)$ denotes the $C^{*}$-algebra of all bounded linear operators on $H$.
For $A \in B (H)$, let $w(A)$ and $\|A\|$ denote the numerical radius and the usual operator norm of $A$, respectively. It is well known that $w(\cdot)$
 defines a norm on $B (H)$, which is equivalent to the usual operator norm $\|\cdot\|$.
 In fact, for every $A\in B(H)$,
\begin{eqnarray}\label{1}
\frac{1}{2}\| A \| \leq w(A) \leq \| A\|.
\end{eqnarray}
An important inequality for $w(A)$ is the power inequality stating that
\begin{eqnarray}\label{2}
w(A^{n})\leq w^{n}(A)
\end{eqnarray}
for each $n\in \mathbb{N}$.
Many authors have investigated several inequalities involving numerical radius inequalities, see e.g. \cite{Bou, drag4, Gol, kit1, She, Zam}.
If $ x,y \in H$ are arbitrary, then the angle between $x$ and $y$ is defined by
 $$\cos \phi _{x,y} =\frac{{\rm Re}\langle x,y\rangle}{\|x\|\|y\|}$$
or by
 $$\cos \psi _{x,y} =\frac{|\langle x,y\rangle |}{\|x\|\|y\|}.$$
The following inequality for angles between two
vectors was obtained by Kre$\breve{i}$n \cite{kre}
\begin{eqnarray}\label{3}
\phi _{x,z}\leq \phi _{x,y}+\phi _{y,z}
\end{eqnarray}
for any nonzero elements $x, y, z \in H.$ By using the representation
\begin{eqnarray}
\psi _{x,y}
= \inf_{\lambda,\mu\in \mathbb{C }- \left\{0\right\} }\phi _{ \lambda x,\mu y} =\inf_{\lambda \in \mathbb{C }- \left\{0\right\} }\phi _{\lambda x, y}=\inf_{\mu \in \mathbb{C }- \left\{0\right\} }\phi _{ x,\mu y}
\end{eqnarray}
and inequality (\ref{3}), he showed that the following triangle inequality is valid
\begin{eqnarray}\label{4}
\psi _{x,y}\leq \psi _{x,z}+\psi _{y,z}
\end{eqnarray}
for any nonzero elements $x,y,z \in H.$

 In section 2  of this paper, we  first introduce some new refinements of numerical radius inequality
 (\ref{1}) by applying the Kre\u in-Lin triangle inequality (\ref{3}) and obtain a reverse of inequality (\ref{2}) in the case when $n=2$.
 In section 3, we obtain some refinements of inequality (\ref{1}) by applying a refinement of the Young inequality.

\section{ Some refinements of inequality (\ref{1}) by Kre\u in-Lin triangle inequality }
In order to achieve our goals, we need the following lemmas.
The first lemma is a simple consequence of the classical Jensen and Young  inequalities.
\begin{lemma}\cite[Lemma 2.1]{moslehian} \label{5}
Let $a,b \geq 0$ and $0 \leq \lambda \leq 1$. Then
$$a^\lambda b^{1-\lambda} \leq \lambda a+(1-\lambda)b\leq [\lambda a^r+(1-\lambda) b^r]^\frac{1}{r},$$
for any $r\geq 1$.
\end{lemma}
The second lemma is a simple consequence of the classical Jensen inequality for convex function $f(t)=t^{r}$, where $r\geq 1.$
\begin{lemma}\label{6}
If $a$ and $b$ are nonnegative real numbers, then
$$(a+b)^{r}\leq 2^{r-1}(a^{r}+b^{r})$$
  for any  $r\geq 1$.
\end{lemma}

\begin{lemma}\cite[Lemma 2.4]{drag3}\label{7}
Suppose that $x,y \in H$  with $\|y\|= 1$.  Then
$$\|x\|^2-|\langle x,y \rangle|^2=\inf_{\lambda \in \mathbb{C}} \| x-\lambda y\|^2.$$
\end{lemma}

 The following lemma is known as a generalized mixed Schwarz inequality.
\begin{lemma}\cite[Lemma 2.3]{moslehian}\label{8}
Let $A \in B(H)$ and $x,y \in H$ be two vectors.
\begin{itemize}
\item[(i)] If
$~~0 \leq \lambda \leq 1$, then
      $$|\langle Ax,y \rangle|^{2}\leq \langle | A|^{2\lambda}x,x \rangle \langle | A^{*}|^{2(1-\lambda)}y,y \rangle.$$
\item[(ii)] If
 $f$ and $g$ are nonnegative continuous functions on $[0,\infty)$ satisfying  $f(t)g(t)=t$, then
    $$|\langle Ax,y \rangle|\leq \|f(| A|)x\| \|g(| A^*|)y\|.$$
\end{itemize}
\end{lemma}
In the next result, we use some ideas of \cite{drag2}.
\begin{theorem}\label{9}
Let $A\in B(H)$ and $f,g$ be nonnegative continuous functions on $[0,\infty)$ satisfying $f(t)g(t)=t$.
 Then for $r \geq 1$
 \begin{eqnarray}\label{55}
 w^{2r}(A)\leq \frac{1}{2^r}\bigg(\| f^2(| A^2|)+g^2(| (A^2)^*|)\|^{r}+2^{r}\inf _{\lambda \in \mathbb{C}}\| A -\lambda I\|^{2r}\bigg).
 \end{eqnarray}
\begin{proof}
By (\ref{4}), we get the  inequality (9) of \cite{drag1} as follows
\begin{eqnarray}\label{10}
\frac{|\langle x,z \rangle|}{\|x\|\|z\|}\frac{|\langle y,z \rangle|}{\|y\|\|z\|}\leq \frac{|\langle x,y \rangle|}{\|x\|\|y\|}+\sqrt{1-\frac{|\langle x,z \rangle|^2}{\|x\|^2\|z\|^2}}\sqrt{1-\frac{|\langle y,z \rangle|^2}{\|y\|^2\|z\|^2}}
\end{eqnarray}
for any $x, y, z \in H\setminus\{0\}$.\\
If we multiply (\ref{10}) by $\|x\| \|y\| \|z\|^2$, then we deduce 
\begin{eqnarray}\label{11}
| \langle x,z \rangle ||\langle y,z \rangle |\leq |\langle x,y \rangle |\|z\|^2+\sqrt{\|x\|^2\|z\|^2- |\langle x,z \rangle|^2}\sqrt{\|y\|^2\|z\|^2-|\langle y,z \rangle|^2}.
\end{eqnarray}
Applying  Lemma \ref{7} for any $x,y,z \in H$ with $\|z\|=1$, we obtain
\begin{eqnarray}\label{12}
 | \langle x,z \rangle||\langle y,z \rangle|\leq |\langle x,y \rangle|+ \inf _{\lambda \in\mathbb{C}}\|x-\lambda z\| \inf_{\mu \in\mathbb{C}} \|y-\mu z\|.
 \end{eqnarray}
 Put
  $x=Az,y = A^{*}z$ in (\ref{12}) to get
\begin{eqnarray}\label{13}
 |\langle Az,z \rangle|^2 & \leq & |\langle A^{2}z,z \rangle|+\inf _{\lambda \in\mathbb{C}}\|Az-\lambda z\| \inf_{\mu \in\mathbb{C}} \|A^{*}z-\mu z\| \nonumber\\ &\leq &
|\langle A^{2}z,z \rangle|+\|Az -\lambda z\|\|A^{*}z-\mu z\|
\end{eqnarray}
for any $z \in H$ with $\|z\|=1$ and $\lambda , \mu \in \mathbb{C}.$\\
 On the other hand, by applying Lemma \ref{8} and the AM-GM inequality, we have
  \begin{eqnarray}\label{14}
 |\langle A^2z,z \rangle| & \leq & \|f(| A^2|)z\| \|g(|(A^2)^{*}|)z\| \nonumber\\ &=& \sqrt{\langle f^2(| A^2|) z,z \rangle\langle g^2(| (A^2)^{*}|)z,z\rangle}\nonumber\\
 &\leq & \frac{1}{2}\langle (f^2(| A^2|)+g^2(| (A^2)^{*}|)) z,z \rangle.
 \end{eqnarray}
 Applying again the AM-GM inequality, we get
  \begin{eqnarray}\label{15}
 \|Az-\lambda z\|\|A^{*}z-\mu z\| \leq \dfrac{\|A z-\lambda z\|^{2}+\|A^*z-\mu z\|^{2}}{2}.
 \end{eqnarray}
 By combining inequalities (\ref{9}), (\ref{14}) and (\ref{15}), we reach
 \begin{eqnarray*}
 |\langle Az,z \rangle|^{2}&\leq & \frac{1}{2}\bigg(\langle( f^{2}(| A^2|)+g^{2}(| (A^{2})^{*}|)) z,z \rangle +\|A z-\lambda z\|^2+\|A^{*} z-\mu z\|^{2}\bigg )\\
 &\leq & \frac{1}{2}\bigg(\langle (f^2(| A^2|)+g^2(| (A^2)^{*}|)) z,z \rangle^{r} +(\|A z-\lambda z\|^{2}+\|A^{*}z -\mu z\|^{2})^{r}\bigg)^{\frac{1}{r}}\\
 &&\hspace{8cm}({\rm by\quad Lemma}~\ref{5})\\
& \leq & \frac{1}{2}\bigg(\langle (f^2(| A^2|)+g^2(| (A^2)^{*}|)) z,z \rangle^{r} +2^{r-1}(\|A z-\lambda z\|^{2r}+\|A^{*}z -\mu z\|^{2r})\bigg)^{\frac{1}{r}}\\ &&\hspace{8cm}({\rm by\quad Lemma} ~\ref{6}).\\
\end{eqnarray*}
Hence
 $$|\langle Az,z \rangle|^{2r} \leq \frac{1}{2^r}\bigg(\langle (f^{2}(| A^{2}|)+g^{2}(| (A^{2})^{*}|)) z,z \rangle^{r} +2^{r-1}(\|A z-\lambda z\|^{2r}+\|A^{*}z-\mu z\|^{2r})\bigg).$$
 By taking the supremum over $z \in H$  with $\| z\| =1 $, we deduce
 $$w^{2r}(A)\leq \frac{1}{2^r}\bigg(\| f^{2}(| A^{2}|)+g^{2}(| (A^{2})^{*}|)\|^{r}+2^{r-1}(\|A-\lambda I\|^{2r}+\|A^{*}-\mu I\|^{2r})\bigg)$$
 for any
 $\lambda,\mu \in \mathbb{C}$.

 Finally, taking the infimum over
  $\lambda,\mu \in \mathbb{C}$ in the inequality above and utilizing
 $$\inf _{\mu \in \mathbb{C}}\|A^{*}-\mu I\| =\inf _{\mu \in \mathbb{C}}\|A-\overline{\mu} I\|=\inf _{\lambda \in \mathbb{C}}\|A-\lambda I\|$$
we obtain the desired result  \eqref{55}.
\end{proof}
\end{theorem}

\begin{remark}\label{16}
In Theorem \ref{9} if we choose
 $r = 1,~~~ f(t)=g(t)=\sqrt{t}$,
we get
$$w^{2}(A)\leq \frac{1}{2}\bigg(\| | A^{2}|+| (A^{2})^{*}|\|+2\inf_{\lambda \in \mathbb{C}} \|A-\lambda I\|^2\bigg).$$
Now, suppose that $s>0$ such that
 $s \leq \sqrt{\|A\|^{2}-\frac{1}{2}\|| A^{2}|+| (A^{2})^{*}|\|}$, if there is $\lambda_{0} \in \mathbb{C}$ in which
 $\|A-\lambda_{0} I\| \leq s$,
 then $w(A)\leq \sqrt{\frac{1}{2}\|| A^{2}|+| (A^{2})^{*}|\|+s^{2}} \leq\|A\|,$
that is an improvement of inequality (\ref{1}) for nonnormal operators.
\end{remark}
Recall that if $A\in M_{2}(\mathbb{R})$, then $\|A\| =\displaystyle{\max_{1\leq i\leq n}\sigma_{i}}$,
 where $\sigma_{i}'s$ are the square root of
eigenvalues of $A^{*}A$, which are called the singular values of $A$,
and $w(A)$ for matrix of the form 
$A=\begin{bmatrix}
a_{1}&b\\
0&a_{2}
\end{bmatrix}$ or 
 $A=\begin{bmatrix}
a_{1}&0\\
b&a_{2}
\end{bmatrix}$ is defined by\\
$$ w(A)=\dfrac{1}{2}|a_{1}+a_{2}|+\dfrac{1}{2}\sqrt{|a_{1}-a_{2}|^{2}+|b|^{2}}, $$
where $a_1, a_2, b\in \mathbb{R}$.
\begin{example}\label{2.7}
By taking
 $A=\begin{bmatrix}
1&\frac{1}{2}\\
0&1
\end{bmatrix}$
and $\lambda_0=\dfrac{1}{2}$
  in Remark \ref{16}, we have $w^{2}(A)\simeq 1.5625,\quad \| A \|^{2}\simeq 3.2822,\quad \|A-\lambda_{0} I\|\simeq 0.5201$
   and $ \frac{1}{2}\| | A^{2}|+| (A^{2})^{*}|\| \simeq 1.5652$.
If $s^{2}\leq \|A\|^{2}-\frac{1}{2}\| |A^{2}|+|(A^{*})^{2}| \|\simeq 1.7170$, then
 $ s\leq1.3103$. Hence our inequality $w(A)\leq \sqrt{\frac{1}{2}\|| A^{2}|+| (A^{2})^{*}|\|+s^{2}} \leq\|A\|$ provides an improvement of inequality (\ref{1}).
\end{example}

\begin{remark}\label{2.8}
 If there exists $\lambda_{0} \in \mathbb{C}$ in which
 $\|A-\lambda_{0} I\| \leq s$, then by putting $\lambda=\mu=\lambda_{0}$ and  by taking the supremum over $z\in H$ with $\| z \|=1$ in (\ref{13}), we deduce
$$w^{2}(A) - w(A^{2})\leq \|A - \lambda_{0} I\|\|A^{*}-\lambda_{0} I\|.$$
Therefore
$$w^{2}(A) - w(A^{2})\leq s^{2}.$$
Now, if $\|A-\lambda_{0} I\| \leq s \leq \sqrt{\|A \|^{2}-w(A^{2})}$, we have $w(A) \leq \sqrt{w(A^{2})+s^{2}}\leq \|A\|$,
 that is an improvement of  inequality (\ref{1}).
\end{remark}
\begin{example}
By taking
$A=\begin{bmatrix}
2&-1\\
0&3
\end{bmatrix}$ and $\lambda_{0}=2.5$ in Remark \ref{2.8},
we have $w(A^{2})\simeq 6.4142,\quad \| A \|^{2}\simeq 10.6054, \quad \|A-\lambda_{0}I\|\simeq 0.955$
for $s\leq\sqrt{\|A\|^{2}-w(A^{2})} \simeq 2.0472$.
Hence our inequality $w(A) \leq \sqrt{w(A^{2})+s^{2}}\leq \|A\|$ provides an improvement of inequality (\ref{1}).
\end{example}

Recall that the vector $x\in H$ is orthogonal to $y \in H$ (denote by $x \perp y$), if $\langle x,y \rangle=0$.
Now, an argument similar to the proof of Theorem \ref{9} with the aid of
Lemma \ref{5} and Lemma \ref{7} gives the following proposition:
\begin{proposition}\label{17}
Let $x, y, z \in H$ with $\|z\| = 1$ and
  $\lambda,\mu \in \mathbb{C},~~~ a, b > 0,$ and $r\geq 1$
   such that
  \begin{eqnarray*}
   \|x-\lambda z\|\leq a,\quad \|y-\mu z\|\leq b.
  \end{eqnarray*}
  Then
  \begin{eqnarray}\label{188}
   (|\langle x,z \rangle||\langle y,z \rangle| - |\langle x,y \rangle|)^{r}\leq \dfrac{a^{2r}+b^{2r}}{2}.
  \end{eqnarray}
  In particular if $x\perp y$, then
  \begin{eqnarray}\label{18}
 ( | \langle x,z \rangle | | \langle y,z \rangle |)^{r} \leq \dfrac{a^{2r}+b^{2r}}{2}
  \end{eqnarray}
  for any $r\geq 1$.
  \begin{proof}
  Since $z$ is a unit vector, from (\ref{11})  we have
  \begin{eqnarray*}
    |\langle x,z \rangle | |\langle y,z \rangle |-|\langle x,y \rangle |
   &\leq & \sqrt{\|x\|^2- | \langle x,z \rangle |^2}\sqrt{\|y\|^{2} -|\langle y,z \rangle |^{2} }\\
 &\leq & \frac{1}{2}(\|x\|^{2}-|\langle x,z \rangle|^{2} +\|y\|^2 -|\langle y,z \rangle|^2)\\
&&\hspace{5cm}(\rm{by~~AM-GM~~ inequality})\\
&=&\frac{1}{2}(\inf_{\lambda \in \mathbb{C}} \| x-\lambda z\|^{2}+ \inf_{\mu\in \mathbb{C}} \| y-\mu z\|^{2})\\
 &&\hspace{5cm}(\rm {by ~ Lemma}~\ref{7}) \\
 &\leq & \frac{1}{2}( \| x-\lambda z\|^{2} + \| y-\mu z\|^{2} \\
 &\leq &\dfrac{a^{2}+b^{2}}{2}\\
&\leq & ( \dfrac{a^{2r}+b^{2r}}{2})^{\frac{1}{r}}~~~~~\hspace{2.8cm} (\rm{by ~ Lemma }~\ref{5}).\\
  \end{eqnarray*}
Hence
  $$(|\langle x,z \rangle| |\langle y,z \rangle|-|\langle x,y \rangle|)^{r} \leq \dfrac{a^{2r}+b^{2r}}{2}.$$
 \end{proof}

 \end{proposition}
\begin{corollary}
Let $A \in B(H)$ and $B$ be a nonzero self-adjoin element in $B(H)$, under assumptions of Proposition \rm \ref{17}, if we choose,
$x = Az$ and $y = Bz$ with $\|z\| = 1$ in (\rm \ref{18}) yields
$$( |\langle Az,z \rangle| |\langle Bz,z \rangle |)^{r} \leq  \dfrac{a^{2r}+b^{2r}}{2}.$$
By taking supremum over $z\in H$ with $\|z\|=1$, we get
$$ w^{r}(A) \leq \dfrac{a^{2r}+b^{2r}}{2} \| B \|^{-r}$$
provided $\|A -\lambda I\|\leq a, ~~~\|B -\mu I\|\leq b$  and for any $r \geq1$ and $a,b >0$.
\end{corollary}
Proposition \ref{17} induces several inequalities as special cases, but here we
only focus on the case $r = 1$, i.e.,
\begin{eqnarray}\label{19}
 |\langle x,z \rangle | | \langle y,z \rangle | \leq \dfrac{a^{2}+b^{2}}{2} +|\langle x,y \rangle |,
\end{eqnarray}
whenever $\|x -\lambda z\|\leq a, ~~~\|y -\mu z\|\leq b$ with $\|z\|=1$ and $\lambda, \mu \in \mathbb{C}$.
\begin{remark}
Suppose that the assumptions of Proposition \ref{17} are still valid.\\
$\bullet$  As an application of inequality (\ref{19}) the following
reverse of inequality (\ref{2}) for $ n=2 $, i.e., an upper bound for
 $w^{2}(A)-w (A^{2})$ can be obtained. In fact,
by choosing  $x = Az$ and $y = A^{*}z$ with $ \|z\| = 1$ and taking supremum over $z \in H$ with $\|z\|=1$, we get
\begin{eqnarray*}
w^{2}(A)-w (A^{2}) \leq \dfrac{a^{2}+b^{2}}{2},
\end{eqnarray*}
provided $\| A -\lambda I\| \leq a,~~~ \| A^{*} -\mu I\| \leq b.$\\
$\bullet$ By choosing  $x = Az$
and $y = A^{-1}z$ with $\|z\| = 1$,  in inequality (\ref{19}) and taking supremum over $z \in H$ with $ \|z\|=1 $, we have
$$ K(A;z)-1\leq \dfrac{a^{2}+b^{2}}{2},$$
provided $\|A -\lambda I\|\leq a,~~~\|A^{-1} -\mu I\|\leq b$, where
$  K(A;z)=\langle Az,z \rangle\langle A^{-1} z,z \rangle $, is the Kantorovich functional.\\
\end{remark}

\section{Some refinements of inequality (\ref{1}) by using Young's inequality}
In this section,  we obtain some refinements of inequality (\ref{1}) by applying refinements of the Young inequality.
 The next lemma is an additive refinement of the scalar Young inequality.
\begin{lemma}\cite[Theorem 2.1]{kit2} \label{20}
If $a,b\geq 0$ and $ 0\leq \lambda \leq 1$, then
\begin{eqnarray*}
a^{\lambda}b^{1-\lambda}+r(\sqrt{a}-\sqrt{b})^{2}\leq \lambda a+(1-\lambda)b,
\end{eqnarray*}
where $r=min \{\lambda,1- \lambda \}$.
\end{lemma}
The main result of this section reads as follows.
\begin{theorem}\label{20}
If $A \in B(H),~~~r = \min\{\lambda,1-\lambda\},$ where $ 0\leq \lambda \leq 1$, then
\begin{eqnarray*}
w(A)\leq \dfrac{1-2r}{2}\| |A|+ |A^{*}|\|+ 2r \| A\|.
\end{eqnarray*}
\begin{proof}Let $x\in H$ be a unit vector. Then we have
\begin{align*}
|\langle Ax,x\rangle | &\leq   \sqrt{\langle |A|x,x\rangle \langle |A^{*}|x,x\rangle}\quad (\rm{by~ Lemma }~\ref{8})\\
&=(\langle |A|x,x\rangle^{1-\lambda}\langle |A^{*}|x,x\rangle^{\lambda})^{\frac{1}{2}}
(\langle|A^{*}|x,x\rangle^{1-\lambda}\langle |A|x,x\rangle^{\lambda})^{\frac{1}{2}}\\
&\leq \frac{1}{2}\bigg(\langle |A|x,x \rangle^{1-\lambda}\langle |A^{*}|x,x\rangle^{\lambda}
+\langle |A^{*}|x,x\rangle^{1-\lambda}\langle |A|x,x\rangle^{\lambda}\bigg)\\
 &\hspace{8cm}(\rm{by~~AM-GM~~inequality})\\
&\leq \frac{1}{2}\bigg((1-\lambda) \langle |A|x,x\rangle+\lambda \langle |A^{*}|x,x\rangle- r(\sqrt{\langle |A|x,x\rangle}-\sqrt{\langle |A^*|x,x\rangle})^{2} \\
&\quad+(1-\lambda) \langle |A^{*}|x,x\rangle+\lambda \langle|A|x,x\rangle- r(\sqrt{\langle |A|x,x\rangle}-\sqrt{\langle |A^*|x,x\rangle})^{2}\bigg)\\
&\hspace{8cm}(\rm{by~ Lemma ~\ref{20}})\\
&= \frac{1}{2}\bigg(\langle(|A|+|A^{*}|)x,x\rangle-2r \langle(|A|+|A^{*}|)x,x\rangle)+4r\sqrt{\langle |A|x,x\rangle\langle |A^{*}|x,x\rangle}\bigg),
\end{align*}
so
$$|\langle Ax,x\rangle |+r\langle (|A|+|A^{*}|)x,x \rangle \leq \frac{1}{2}(\langle(|A|+|A^{*}|)x,x\rangle+ 4r\sqrt{\langle |A|x,x\rangle\langle |A^{*}|x,x\rangle}).$$
By taking supremum over $x\in H$ with $\|x\|=1$, we deduce
\begin{eqnarray*}
w(A)\leq \frac{1-2r}{2}\| |A|+|A^{*}|\|+2r \|A\|,
\end{eqnarray*}
which is an improvement of inequality (\ref{1}).
\end{proof}
\end{theorem}

\begin{example}
Let $A=\begin{bmatrix}
1&1\\
0&2
\end{bmatrix}$ be as in Theorem \ref{20} and $r=0.1$. Then by straightforward computation, we get
$w(A)\simeq2.2071, \quad \| A \|\simeq2.2882$ and $ \dfrac{1}{2}\||A|+|A^{*}|\|\simeq2.2518 $. Hence
\begin{eqnarray*}
w(A)\leq \frac{1-2r}{2}\| |A|+|A^{*}|\|+2r \|A\|\leq \|A\|,
\end{eqnarray*}
 provides an improvement of inequality (\ref{1}). In fact, $2.2071\leq 2.2590 \leq 2.2882$.
\end{example}
The following lemma is a multiplicative refinement of the Young inequality with the Kantorovich constant.
\begin{lemma}\cite[Corollary 3]{Fuji}\label{21}
Let $a,b>0$. Then
\begin{eqnarray*}
(1-\lambda) a+\lambda b\geq k(h,2)^{r} a^{1-\lambda}b^{\lambda},
\end{eqnarray*}
where $0 \leq \lambda \leq 1,~~~ r=\min \{\lambda, 1-\lambda\},~~~ h=\frac{b}{a}$  such that $K(h,2)=\frac{(h+1)^{2}}{4h}$ for $h>0$, which has properties
 $K(h,2)=K(\frac{1}{h},2)\geq 1(h>0)$ and $K(h,2)$ is increasing on $[1,\infty)$ and is decreasing on $(0,1)$.
\end{lemma}

In \cite{kit3}, Kittaneh obtained the inequality
\begin{eqnarray}\label{22}
w(A) \leq \dfrac{1}{2}\| |A|+|A^*|\|.
\end{eqnarray}
In the following theorem, we improve  inequality (\ref{22}) for hyponormal operators.
Before proceeding recall that the operator $A\in B(H)$ is said to be hyponormal if $A^{*}A-AA^{*}\geq 0$.

\begin{theorem}
If $A\in B(H)$ is hyponormal, $r=\min \{\lambda,1-\lambda\},$ where $0 \leq \lambda \leq 1 $, then
\begin{eqnarray*}
w(A)\leq \dfrac{1}{\displaystyle{\inf_{\| x \|=1}\zeta(x)}}\dfrac{\| |A|+|A^{*}|\|}{2},
\end{eqnarray*}
where $\zeta(x)=K(\frac{\langle |A|x,x \rangle}{\langle |A^{*}|x,x \rangle},2)^{r}$ is a refinement of inequality (\ref{1}).
\begin{proof}Let $x\in H$ be a unit vector.
\begin{align*}
|\langle Ax,x \rangle | &\leq \sqrt{\langle |A|x,x \rangle \langle |A^{*}|x,x \rangle}\quad \quad(by~ Lemma ~\ref{8})\\
&=(\langle |A^{*}|x,x \rangle^{1-\lambda} \langle |A|x,x \rangle^{\lambda})^{\frac{1}{2}}(\langle |A|x,x \rangle^{1-\lambda} \langle |A^{*}|x,x \rangle^{\lambda})^{\frac{1}{2}}\\
&\leq  \dfrac{1}{2}\bigg((\langle |A^{*}|x,x \rangle^{1-\lambda}\langle |A|x,x \rangle^{\lambda})+(\langle |A|x,x \rangle^{1-\lambda} \langle |A^{*}|x,x \rangle^{\lambda}\bigg)\\
&\hspace{5cm}(\rm{by~~AM-GM~~inequality})\\
&\leq \dfrac{1}{2}\bigg(\dfrac{1}{K(\frac{\langle |A|x,x \rangle}{\langle |A^{*}|x,x \rangle},2)^{r}}((1-\lambda)\langle |A^{*}|x,x \rangle +\lambda \langle |A|x,x \rangle)\\
&\quad+(\dfrac{1}{K(\frac{\langle |A^{*}|x,x \rangle}{\langle |A|x,x \rangle},2)^{r}}((1-\lambda)\langle |A|x,x \rangle +\lambda \langle |A^{*}|x,x \rangle\bigg)\\
 &\hspace{5cm}(\rm{by~Lemma} ~\ref{21})\\
&= \dfrac{1}{2}\bigg(\frac{1}{K(\frac{\langle |A|x,x \rangle}{\langle |A^{*}|x,x\rangle},2)^{r}}(\langle |A^{*}|x,x\rangle + \langle |A|x,x\rangle)\bigg).
\end{align*}
Taking supremum over $x\in H$ with $\|x\|=1$, we have
\begin{eqnarray*}
w(A)\leq \dfrac{1}{\displaystyle{\inf_{\| x \|=1}\zeta(x)}}\dfrac{\| |A|+|A^{*}|\|}{2},
\end{eqnarray*}
where $\zeta(x)=K(\frac{\langle |A|x,x \rangle}{\langle |A^{*}|x,x \rangle},2)^{r}$.\\
Note that $2\langle |A|x,x \rangle\langle |A^{*}|x,x \rangle\leq \langle |A|x,x \rangle^{2}+\langle |A^{*}|x,x \rangle^{2}$,
so
\begin{eqnarray*}
 (\langle |A|x,x \rangle+\langle |A^{*}|x,x \rangle)^{2}\geq 4 \langle |A|x,x \rangle\langle |A^{*}|x,x \rangle.
\end{eqnarray*}
Hence
\begin{eqnarray*}
\dfrac{(\langle |A|x,x \rangle+\langle |A^{*}|x,x \rangle)^{2}}{4\langle |A|x,x \rangle\langle |A^{*}|x,x \rangle}\geq 1.
\end{eqnarray*}
Therefore $K\left(\frac{\langle |A|x,x \rangle}{\langle |A^{*}|x,x \rangle},2\right)\geq 1$.
\end{proof}
\end{theorem}
\bibliographystyle{amsplain}

\end{document}